\documentclass{article}

\usepackage{amsmath,amssymb,amsthm}

\newtheorem{theorem}{Theorem}[section]
\newtheorem{definition}[theorem]{Definition}
\newtheorem{proposition}[theorem]{Proposition}
\newtheorem{corollary}[theorem]{Corollary}
\newtheorem{lemma}[theorem]{Lemma}
\newtheorem{fact}[theorem]{Remark}
\newtheorem{exemplu}[theorem]{Example}

\newcommand{\bdfn}{\begin{definition}}
\newcommand{\edfn}{\end{definition}}
\newcommand{\bthm}{\begin{theorem}}
\newcommand{\ethm}{\end{theorem}}
\newcommand{\bprop}{\begin{proposition}}
\newcommand{\eprop}{\end{proposition}}
\newcommand{\bcor}{\begin{corollary}}
\newcommand{\ecor}{\end{corollary}}
\newcommand{\blem}{\begin{lemma}}
\newcommand{\elem}{\end{lemma}}
\newcommand{\bfact}{\begin{fact}}
\newcommand{\efact}{\end{fact}}
\newcommand{\bex}{\begin{exemplu}\begin{rm}}
\newcommand{\eex}{\end{rm}\end{exemplu}}

\def\R{{\mathbb R}}
\def\N{{\mathbb N}}

\newcommand{\lambdaxy}{(1-\lambda)x\oplus\lambda y}

\newcommand{\midxy}{\frac12x\oplus\frac12y}

\newcommand{\eps}{\varepsilon}

\newcommand{\be}{\begin{enumerate}}
\newcommand{\ee}{\end{enumerate}}
\newcommand{\bt}{\begin{tabular}}
\newcommand{\et}{\end{tabular}}
\newcommand{\beq}{\begin{equation}}
\newcommand{\eeq}{\end{equation}}
\newcommand{\ba}{\begin{array}} 
\newcommand{\ea}{\end{array}}
\newcommand {\bea} {\begin{eqnarray}}
\newcommand {\eea} {\end {eqnarray}}
\newcommand {\bua} {\begin{eqnarray*}}
\newcommand {\eua} {\end {eqnarray*}}
\newcommand{\se}{\subseteq}
\newcommand{\ds}{\displaystyle}

\begin{document}

\title{A quadratic rate of asymptotic regularity for CAT(0)-spaces}

\author{Lauren\c tiu Leu\c stean\\[0.2cm] 
\footnotesize Department of Mathematics, Darmstadt University of Technology,\\
\footnotesize Schlossgartenstrasse 7, 64289 Darmstadt, Germany\\[0.1cm]
\footnotesize and\\
\footnotesize Institute of Mathematics "Simion Stoilow'' of the Romanian Academy, \\
\footnotesize Calea Grivi\c tei 21, P.O. Box 1-462, Bucharest, Romania\\[0.1cm]
\footnotesize E-mail: leustean@mathematik.tu-darmstadt.de}
\date{}
\maketitle

\begin{abstract}
In this paper we obtain a quadratic bound on the rate of asymptotic regularity for the Krasnoselski-Mann iterations of nonexpansive mappings in CAT(0)-spaces, whereas previous results guarantee only exponential bounds.
The method we use is to extend to the more general setting of uniformly convex hyperbolic spaces a quantitative version of a strengthening of Groetsch's theorem obtained by Kohlenbach  using methods from mathematical logic (so-called ``proof mining''). 
\end{abstract}

\begin{tabular}{ll}
\footnotesize\noindent {\it Keywords}: &\footnotesize Proof mining, metric fixed point theory, nonexpansive functions, \\
&\footnotesize asymptotic regularity, CAT(0)-spaces, hyperbolic spaces\\[0.1cm]
\noindent {\it MSC:\ } & \footnotesize 47H10,  47H09, 03F10
\end{tabular}

\section{Introduction}

In this paper we present another case study in the general project of {\em proof mining} in functional analysis, developed by Kohlenbach (see \cite{Kohlenbach-book} for details). By "proof mining" we mean the logical analysis of mathematical proofs with the aim of extracting new numerically relevant information hidden in the proofs.

Thus, we obtain a quadratic bound on the rate of asymptotic regularity for the Krasnoselski-Mann iterations of nonexpansive self-mappings of nonempty convex, bounded sets $C$ in CAT(0)-spaces in the sense of Gromov (see \cite{Bridson/Haefliger} for a detailed treatment). Moreover, the bound we get is uniform in the sense  that does not depend on the nonexpansive mapping $T$ or on  the starting point $x\in C$, and depends on $C$ only via its diameter.

The method we use to get this bound is to find explicit uniform bounds on  the rate of asymptotic regularity in the general setting of uniformly convex hyperbolic spaces, and then to specialize them to CAT(0)-spaces. 

The notion of nonexpansive mapping can be introduced in the very general setting of metric spaces. Thus, if $(X,\rho)$ is a metric space, and $C\subseteq X$ a nonempty subset, then a mapping $T:C\to C$ is called {\em nonexpansive} if for all $x,y\in C$,
\[\rho(Tx, Ty)\le \rho(x,y).\]

Different notions of "hyperbolic space'' \cite{Kirk(82),Goebel/Kirk(83),Goebel/Reich(84),Reich/Shafrir(90)} can be found in the literature. We work in the setting of hyperbolic spaces as introduced by Kohlenbach \cite{Kohlenbach(05a)}, which are slightly more restrictive than the spaces of hyperbolic type \cite{Goebel/Kirk(83)} by (W4), but more general then the concept of hyperbolic space from \cite{Reich/Shafrir(90)}. See \cite{Kohlenbach(05a),Kohlenbach/Leustean(05)} for detailed discussion of this and related notions.

A {\em  hyperbolic space}
is a triple $(X,\rho,W)$ where
$(X,\rho)$ is a metric space and $W:X\times X\times [0,1]\to X$ is such that 
\begin{eqnarray*}
(W1) & \rho(z,W(x,y,\lambda))\le (1-\lambda)\rho(z,x)+\lambda \rho(z,y),\\
(W2) & \rho(W(x,y,\lambda),W(x,y,\tilde{\lambda}))=|\lambda-\tilde{\lambda}|\cdot 
\rho(x,y),\\
(W3) & W(x,y,\lambda)=W(y,x,1-\lambda),\\
(W4) & \,\,\,\rho(W(x,z,\lambda),W(y,w,\lambda)) \le (1-\lambda)\rho(x,y)+\lambda
\rho(z,w).
\end {eqnarray*}

\noindent If $x,y\in X$, and $\lambda\in[0,1]$ then we use the notation $(1-\lambda)x\oplus \lambda y$ for $W(x,y,\lambda)$. 

We remark that any normed space $(X,\|\cdot\|)$ is a hyperbolic space, with $(1-\lambda)x\oplus\lambda y:=(1-\lambda)x+\lambda y$.

The notion of uniformly convex hyperbolic space $(X,\rho,W)$ with a modulus of uniform convexity $\eta$ is defined in Section \ref{UC-hyp} following the normed space case.

A very important class of hyperbolic spaces are the CAT(0)-spaces. Thus, a hyperbolic space is a CAT(0)-space if and only if it satisfies the so-called CN-inequality  of Bruhat-Tits \cite{Bruhat/Tits}: for all $x,y,z\in X$,
\begin{equation}
\quad \rho\left(z,\midxy\right)^2\le \frac12\rho(z,x)^2+\frac12\rho(z,y)^2-\frac14\rho(x,y)^2.\label{CN}
\end{equation}
(see \cite[p.163]{Bridson/Haefliger}, and \cite[p. 98]{Kohlenbach(05a)} for details).

Moreover, it will turn out that CAT(0)-spaces are uniformly convex hyperbolic spaces with a quadratic modulus of convexity.

In the sequel, $(X,\rho,W)$ is a hyperbolic space, $C\se X$ a nonempty convex subset of $X$, and $T:C\to C$ a nonexpansive mapping. 

As in the case of normed spaces \cite{Mann(53),Krasnoselski(55)},  we can define the {\em Krasnoselski-Mann iteration} starting from $x\in C$ by: 
\begin{equation}
x_0:=x, \quad x_{n+1}:=(1-\lambda_n)x_n \oplus\lambda_n Tx_n, \label{KM-lambda-n-def-hyp}\end{equation}
where $(\lambda_n)$ is a sequence in $(0,1)$.

Asymptotic regularity  was already implicit in \cite{Krasnoselski(55),Schaefer(57),Edelstein(66)}), but only in 1966 Browder and Petryshyn \cite{Browder/Petryshyn(66)} defined it for normed spaces $(X,\|\cdot\|)$. In our setting, the mapping $T:C\to C$ is called {\em asymptotic regular}  if  for all $x\in C$,
\[\ds\lim_{n\to\infty}\rho(T^n(x),T^{n+1}(x))=0.\]
  
For constant $\lambda_n=\lambda\in(0,1)$, the fact that $\ds\lim_{n\to\infty}\rho(x_n,Tx_n)=0$ for all $x\in C$ is equivalent to the asymptotic regularity of the averaged mapping 
\[T_\lambda:=(1-\lambda)I\oplus\lambda T.\] 

Therefore, for  general $\lambda_n\in(0,1)$, following \cite{Borwein/Reich/Shafrir(92)}, we say that the nonexpansive mapping $T$ is {\em $\lambda_n$-asymptotic regular} if  for all $x\in C$, 
\[\lim_{n\to\infty}\rho(x_n,Tx_n)=0.\]

The most general assumptions on the sequence $(\lambda_n)$ for which asymptotic regularity has been proved for arbitrary normed spaces  and bounded sets $C$ are the following: 
\begin{eqnarray}
\sum_{n=0}^\infty \lambda_n=+\infty, \text{~and}\label{lambda-n-div}\\
\limsup_{n\to \infty}\lambda_n<1. \label{lambda-n-limsup}
\end{eqnarray}

Thus, Ishikawa proved in his seminal paper \cite{Ishikawa(76)} one of the most important results in the fixed point theory of nonexpansive mappings.

\begin{theorem}
Assume that $(\lambda_n)$ satisfies (\ref{lambda-n-div}), (\ref{lambda-n-limsup}), and that $C$ is bounded. Then $T$ is $\lambda_n$-asymptotic regular.
\end{theorem}
Independently, Edelstein/O'Brien \cite{Edelstein/OBrien(78)} also proved the asymptotic regularity for constant $\lambda_n=\lambda\in(0,1)$, and noted that it is uniform in $x$.

By a logical analysis of the proof of a theorem due to Borwein/Reich/Shafrir \cite{Borwein/Reich/Shafrir(92)} (which generalizes Ishikawa's result to unbounded $C$), Kohlenbach \cite{Kohlenbach(01)} obtained for the first time explicit bounds on the asymptotic regularity when general sequences $(\lambda_n)$ satisfying (\ref{lambda-n-div}), (\ref{lambda-n-limsup}) are considered.

Subsequently, Kohlenbach and the author \cite{Kohlenbach/Leustean(03)} extended these results to the very general setting of hyperbolic spaces (and even to the more general class of directional nonexpansive mappings as introduced in \cite{Kirk(00)}). 

The following result (which is a corollary of the main theorem in \cite{Kohlenbach/Leustean(03)}) is proved there for hyperbolic spaces in the sense of \cite{Reich/Shafrir(90)}, but the proof goes through for the setting used in the present paper.

\begin{theorem}\label{bounds-Ishikawa}
Let $(X,\rho, W)$ be a hyperbolic space, $C\subseteq X$ a nonempty
convex bounded subset with diameter $d_C$, and $T:C\rightarrow C$ a nonexpansive mapping. \\
Assume that $K\in\N, K\geq 2$ and $(\lambda_n)$ is a
sequence in $\left[\ds\frac{1}{K}, 1-\frac{1}{K}\right]$.\\
Then $T$ is $\lambda_n$-asymptotic  regular, and the following holds for all $x\in C$:
\[ \forall \varepsilon >0\,\forall n\ge h(\varepsilon,d,K)\,
\big(\rho(x_n,Tx_n) <\varepsilon\big), \]
where
\[\begin{array}{l} 
\ds h(\varepsilon,d,K):=K\cdot M\cdot\lceil 2d\cdot\exp(K(M+1))\rceil,
\mbox{ with} \\
\ds d\in\R, \  d\ge d_C,\ \mbox{and} \\
M\in \N, \ M\ge\ds\frac{1+2d}{\varepsilon}. 
\end{array}\]
\end{theorem}

For normed spaces and the special case of constant $\lambda_n=\lambda\in (0,1)$ the exponential bound in the above theorem is not optimal. In this case, a uniform and optimal quadratic bound was obtained by Baillon/Bruck \cite{Baillon/Bruck(96)} using an extremely complicated computer aided proof, and only for $\lambda_n=1/2$ a classical proof of a  result of this type was given \cite{Bruck(96)}.
However, the questions whether the methods of proof used by them hold for  non-constant sequences $\lambda_n$ or for hyperbolic spaces are left as open problems in \cite{Baillon/Bruck(96)}, and as far as we know they received no positive answer until now. Hence, the bound from Theorem \ref{bounds-Ishikawa} is the only effective bound known at all for non-constant sequences $\lambda_n$ (even for normed spaces).

Our result guarantees only an exponential bound for the asymptotic regularity in the case of CAT(0)-spaces, and as we have already remarked, it seems that Baillon/Bruck's approach does not extend to the more general setting of hyperbolic spaces.

In this paper we show that we can still get a quadratic rate of asymptotic regularity for CAT(0)-spaces, but following a completely different approach, inspired by the results on asymptotic regularity obtained before Ishikawa, and Edelstein/O'Brien theorems, in the setting of uniformly convex normed spaces. 

More specifically, our point of departure is the following result, proved by Groetsch \cite{Groetsch(72)} (see also \cite{Reich(79)}):  

\begin{theorem}\label{Groetsch-thm}
Let $(X,\|\cdot\|)$ be a uniformly convex normed space, $C\se X$ a nonempty convex subset, and $T:C\to C$ a nonexpansive mapping such that $T$ has at least one fixed point. \\
Assume that $(\lambda_n)$ satisfies the following condition
\beq
\sum_{k=0}^\infty \lambda_k(1-\lambda_k)=\infty.\label{hyp-lambda-n-Groetsch}
\eeq 
Then $T$ is $\lambda_n$-asymptotic regular.
\end{theorem}

The assumption  (\ref{hyp-lambda-n-Groetsch}) is equivalent with the existence of a witness $\theta:\N\to\N$ such that for all $n\in \N$,
\beq
\sum_{k=0}^{\theta(n)} \lambda_k(1-\lambda_k)\geq n.\label{def-theta}
\eeq

By ``proof mining'', Kohlenbach \cite{Kohlenbach(03)} obtained a quantitative version of a strengthening of Groetsch's theorem which only assumes the existence of approximate fixed points in some neighborhood of $x$, and generalizing previous results obtained by Kirk/Martinez-Yanez \cite{Kirk/Martinez(90)} for constant $\lambda_n=\lambda\in(0,1)$. 

In  this paper, we  extend Kohlenbach's results  to the more general setting of uniformly convex hyperbolic spaces. 
In this way, for bounded $C$, we get $\lambda_n$-asymptotic regularity for general $\lambda_n$ satisfying (\ref{hyp-lambda-n-Groetsch}), and we also obtain explicit bounds on the rate of asymptotic regularity, which are uniform in $x,T$, and depend on the uniformly convex hyperbolic space $(X,\rho,W)$ only via a modulus $\eta$, on $C$ only weakly via its diameter $d_C$, and on $\lambda_n$ only via the witness $\theta$.

The most important consequence of our results is that for CAT(0)-spaces, which are uniformly convex hyperbolic spaces with a nice quadratic modulus $\eta$, we obtain a quadratic rate of asymptotic regularity (see Corollary \ref{CAT0-constant-lambda}).

The following table presents a general picture of the cases where effective bounds for asymptotic regularity were obtained (UC means uniformly convex).\\[0.4cm]
\begin{tabular}{|l|l|l|}
\hline 
& $\lambda_n=\lambda$ & non-constant $\lambda_n$ \\[0.2cm]
\hline
Hilbert spaces & quadratic: \cite{Browder/Petryshyn(67)} & $\ds\theta\left(\frac 1{\varepsilon^2}\right)$: \cite{Kohlenbach(03)}\\
\hline
$\l_p$ spaces, $2\leq p<\infty$ &  quadratic: \cite{Kohlenbach(03),Kirk/Martinez(90)} & $\ds\theta\left(\frac 1{\varepsilon^p}\right)$: \cite{Kohlenbach(03)}\\
\hline
UC normed spaces &\cite{Kohlenbach(03),Kirk/Martinez(90)} & \cite{Kohlenbach(03)}\\
\hline
normed spaces & quadratic: \cite{Baillon/Bruck(96)} & \cite{Kohlenbach(01)}\\
\hline
CAT(0)-spaces & quadratic: present paper & $\ds\theta\left(\frac 1{\varepsilon^2}\right)$:  present paper\\
\hline
UC hyperbolic spaces & present paper & present paper\\
\hline
hyperbolic spaces & \cite{Kohlenbach/Leustean(03)} & \cite{Kohlenbach/Leustean(03)} \\
\hline
\end{tabular}

\section{Uniformly convex hyperbolic spaces}\label{UC-hyp}

In the following, $(X,\rho,W)$ is a hyperbolic space. 

\begin{lemma}
Let $x,y\in X$, and $\lambda\in[0,1]$. Then 
\begin{equation}
\rho(x,\lambdaxy)=\lambda\rho(x,y),\text{~and~} \rho(y,\lambdaxy)=(1-\lambda)\rho(x,y).\label{prop-xylambda}
\end{equation}
\end{lemma}
\begin{proof}
Applying (W1) with $z:=x,y$, we get that 
\begin{equation}
\rho(x,\lambdaxy)\le \lambda\rho(x,y), \quad \rho(y,\lambdaxy)\le (1-\lambda)\rho(x,y).\label{eq1}
\end{equation} 
It follows that
\begin{eqnarray*}
\rho(x,y)&\le &\rho(x,\lambdaxy)+\rho(y,\lambdaxy)\\
&\le& \lambda\rho(x,y)+(1-\lambda)\rho(x,y)
=\rho(x,y),
\end {eqnarray*}
hence we must have equality in (\ref{eq1}).
\qed\end{proof}

Following \cite{Takahashi(70)},  $(X,\rho,W)$ is called {\em strictly convex} if for any $x,y\in X,$ and $\lambda\in[0,1]$ there exists a unique element $z\in X$ such that 
\[\rho(z,x)=\lambda\rho(x,y), \quad \text{and~}\rho(z,y)=(1-\lambda)\rho(x,y).\]

We define uniform convexity following \cite[p.105]{Goebel/Reich(84)}.

\begin{definition}
The hyperbolic space $(X,\rho,W)$ is called {\em uniformly convex} if for 
any $r>0$, and $\varepsilon\in(0,2]$ there exists a $\delta\in(0,1]$ such that 
for all $a,x,y\in X$,
\begin{eqnarray}
\left.\begin{array}{l}
\rho(x,a)\le r\\
\rho(y,a)\le r\\
\rho(x,y)\ge\varepsilon r
\end{array}
\right\}
& \quad \Rightarrow & \quad \rho\left(\frac12x\oplus\frac12y,a\right)\le (1-\delta)r. \label{uc-def}
\end{eqnarray}
A mapping $\eta:(0,\infty)\times(0,2]\rightarrow (0,1]$ providing such a
$\delta:=\eta(r,\varepsilon)$ for given $r>0$ and $\varepsilon\in(0,2]$ is called a {\em modulus of uniform
convexity}.
\end{definition}

The proofs of the following facts about uniform convex hyperbolic spaces are similar to the ones of the corresponding results for uniformly convex Banach spaces. In order to make the paper self-contained, we include them. 

\begin{proposition}\label{uc=>sc}
Any uniformly convex hyperbolic space is strictly convex.
\end{proposition}
\begin{proof}
Let $(X,\rho,W)$ be uniformly convex with modulus of uniform convexity $\eta$. We proof that $(X,\rho,W)$ is strictly convex by contradiction. Thus, assume that there are $x,y\in X$, and $\lambda\in[0,1]$ such that there exist two distinct points $z,w\in X$ with 
\[\rho(z,x)=\rho(w,x)=\lambda\rho(x,y), \quad \rho(z,y)=\rho(w,y)=(1-\lambda)\rho(x,y).\]
It follows immediately that $x\ne y$ and $\lambda\in (0,1)$. Let $\displaystyle r_1:=\lambda\rho(x,y)>0, r_2:=(1-\lambda)\rho(x,y)>0, \varepsilon_1:=\frac{\rho(z,w)}{r_1},$ and $\displaystyle\varepsilon_2:=\frac{\rho(z,w)}{r_2}$. It is easy to see that $\varepsilon_1,\varepsilon_2\in(0,2]$, so we can apply twice (\ref{uc-def}) to get 
\[\rho\left(\frac12z\oplus\frac12w,x\right)\le (1-\eta(r_1,\varepsilon_1))r_1,\quad \rho\left(\frac12z\oplus\frac12w,y\right)\le (1-\eta(r_2,\varepsilon_2))r_2.\]
Since $x\ne y$, we must have $\eta(r_1,\varepsilon_1)<1$ or $\eta(r_2,\varepsilon_2)<1$.
It follows that
\begin{eqnarray*}
\rho(x,y)&\le & \rho\left(\frac12z\oplus\frac12w,x\right)+\rho\left(\frac12z\oplus\frac12w,y\right)\\&\le&  (1-\eta(r_1,\varepsilon_1))r_1+(1-\eta(r_2,\varepsilon_2))r_2<r_1+r_2=\rho(x,y),
\end {eqnarray*}
that is a contradiction.
\qed\end{proof}

\begin{proposition}
Let $\eta:(0,\infty)\times(0,2]\to(0,1]$. 
The following are equivalent:
\begin{enumerate}
\item $(X,\rho,W)$ is uniformly convex with modulus of uniform convexity $\eta$;
\item for any $r>0,\varepsilon\in(0,2], \lambda\in[0,1]$, and $a,x,y\in X$, 
\begin{eqnarray}
\left.\begin{array}{l}
\rho(x,a)\le r\\
\rho(y,a)\le r\\
\rho(x,y)\ge\varepsilon r
\end{array}
\right\}
& \quad \Rightarrow & \quad\rho(\lambdaxy,a)\le (1-\gamma(r,\varepsilon,\lambda))r,\label{uc-lambda}
\end {eqnarray}
where 
\[\gamma(r,\varepsilon,\lambda)=
\begin{cases}
2\lambda\eta(r,\varepsilon), & \text{if~} \lambda\le \frac12\\
2(1-\lambda)\eta(r,\varepsilon), & \text{otherwise.}
\end{cases}\]
\end{enumerate}
\end{proposition}
\begin{proof}
$(ii)\Rightarrow(i)$ Just take $\displaystyle\lambda:=\frac12$.\\
$(i)\Rightarrow(ii)$ Assume first that $\displaystyle\lambda\le \frac12$. Applying (\ref{prop-xylambda}) and (W2), we get that
\begin{eqnarray*}
\rho(\lambdaxy,x)&=&\lambda\rho(x,y)=(2\lambda)\rho\left(\midxy,x\right),\\
\rho\left(\lambdaxy,\midxy\right)&=&\left|\lambda-\frac12\right|\rho(x,y)=(1-2\lambda)\rho\left(\midxy,x\right).
\end {eqnarray*}
Since $X$ is also strictly convex, we must have 
\[\lambdaxy=(1-2\lambda)x\oplus(2\lambda)\left(\midxy\right).\]
 It follows then
\begin{eqnarray*}
\rho(\lambdaxy,a)&=&\rho\left((1-2\lambda)x\oplus(2\lambda)\left(\midxy \right),a\right)\\
&\le &(1-2\lambda)\rho(x,a)+2\lambda\rho\left(\midxy,a \right)\\
&\le& (1-2\lambda)r+2\lambda(1-\eta(r,\varepsilon))r =(1-2\lambda\eta(r,\varepsilon))r. 
\end {eqnarray*}
If $\displaystyle\lambda> \frac12$, then $\displaystyle 1-\lambda<\frac12$, so
\[\rho(\lambda y\oplus(1-\lambda)x,a)\le(1-2(1-\lambda)\eta(r,\varepsilon))r.\] 
Use now (W3) to get the conclusion.
\qed\end{proof}

\begin{lemma}
Let $(X,\rho,W)$ be a uniformly convex hyperbolic space with modulus of uniform convexity $\eta$.
For any $r>0,\varepsilon\in(0,2], \lambda\in[0,1]$, and $a,x,y\in X$, 
\begin{eqnarray}
\left.\begin{array}{l}
\rho(x,a)\le r\\
\rho(y,a)\le r\\
\rho(x,y)\ge\varepsilon r
\end{array}
\right\}
&\quad \Rightarrow & \quad\rho(\lambdaxy,a)\le  (1-2\lambda(1-\lambda)\eta(r,\varepsilon))r. \label{uc-ineq-Groetsch}
\end {eqnarray}
\end{lemma}
\begin{proof}
Apply (\ref{uc-lambda}) and the fact that $2\lambda, 2(1-\lambda)\ge 2\lambda(1-\lambda)$.
\qed\end{proof}

\begin{proposition}\label{CAT0-uc}
Assume that $X$ is a CAT(0)-space. Then $X$ is uniformly convex, and 
\begin{equation}
\eta(r,\varepsilon):=\frac{\varepsilon^2}{8}\label{mod-conv-CAT0}
\end{equation}
is a modulus of uniform convexity such that 
\begin{enumerate}
\item $\eta$ decreases with  $r$(for a fixed $\varepsilon$).
\item $\eta(r,\varepsilon)$ can be written as  $\eta(r,\varepsilon)=\varepsilon\cdot\tilde{\eta}(r,\varepsilon)$,and $\tilde{\eta}$ increases with $\varepsilon$ (for a fixed $r$).
\end{enumerate}
\end{proposition}
\begin{proof}
Let $r>0$, $\varepsilon\in(0,2]$,  and $a,x,y\in X$ be such that $\rho(x,a)\le r,\rho(y,a)\le r ,\rho(x,y)\ge\varepsilon r$. 
Applying (\ref{CN}) we get that 
\begin{eqnarray*}
\rho\left(\midxy,a\right)\!\!\!\!\!&\le &\!\!\!\!\sqrt{\frac12\rho(x,a)^2+\frac12\rho(y,a)^2-\frac14\rho(x,y)^2}\le \sqrt{\frac12 r^2+\frac12 r^2-\frac14\varepsilon^2r^2}\\
&=& \sqrt{1-\frac{\varepsilon^2}4}\cdot r\le \left(1-\frac{\varepsilon^2}8\right)r.
\end {eqnarray*}
Hence, $X$ is uniformly convex, and  $\eta(r,\varepsilon)$  defined by (\ref{mod-conv-CAT0}) is a modulus of uniform convexity for $X$. (i) and (ii) follow immediately.
\qed\end{proof}

\section{Technical results}

In the sequel, $(X,\rho,W)$ is a uniformly convex hyperbolic space with modulus of uniform convexity $\eta$, $C\subseteq X$ a nonempty convex subset, and $T:C\to C$ nonexpansive.

\begin{lemma}\label{useful-KM-x-y}
Let $x,y\in X$, and $(x_n)$ be the Krasnoselski-Mann iterating starting with $x$. Then 
\begin{enumerate}
\item $(\rho(x_n,Tx_n))$ is non-increasing;
\item for any $n\in\N$
\begin{eqnarray}
\rho(x_{n+1},y)&\le& \rho(x_n,y)+\lambda_n\rho(y,Ty); \label{KM-useful-1}\end {eqnarray}
\item
\begin{eqnarray}\rho(x_n,y)&\le& \rho(x,y)+\left(\sum_{i=0}^{n-1}\lambda_i\right)\rho(y,Ty). \label{KM-useful-2}
\end {eqnarray}
\end{enumerate}
\end{lemma}
\begin{proof}
\begin{enumerate}
\item See \cite[Proposition 3.4]{Kohlenbach/Leustean(03)}, whose proof immediately generalizes to the notion of hyperbolic space used in this paper. 
\item
\begin{eqnarray*}
\rho(x_{n+1},y)&\le & (1-\lambda_n)\rho(x_n,y)+\lambda_n\rho(Tx_n,y)\\
&\le &(1-\lambda_n)\rho(x_n,y)+\lambda_n\rho(Tx_n,Ty)+\lambda_n\rho(Ty,y)\\
&\le &\rho(x_n,y)+\lambda_n\rho(Ty,y).
\end {eqnarray*}
\item (\ref{KM-useful-2})  follows from (\ref{KM-useful-1}) by an easy induction.
\end{enumerate}
\qed\end{proof}

For any $x,y\in X, n\in\N$, let us denote 
\begin{equation}
\alpha(x,n,y):=\rho(x_n,y)+\rho(y,Ty). \label{alpha-def}
\end{equation}

\begin{lemma}{\bf (Main technical lemma)}\label{main-lemma}\\
Assume that $\eta$ decreases with $r$ (for a fixed $\varepsilon$).\\ 
Let $x,y\in X$, $n\in\N$, and  $\gamma,\beta,\tilde{\beta},a>0$ be such that 
\begin{equation}
\gamma\le \alpha(x,n,y)\le \beta,\tilde{\beta}, \quad  \text{and~~} a\le\rho(x_n,Tx_n).\label{hyp-main-lemma}
\end{equation}
Then the following inequality is satisfied:
\begin{eqnarray}
\rho(x_{n+1},y)&\le& \rho(x_n,y)+\rho(y,Ty)- 2\gamma\lambda_n(1-\lambda_n)\eta\left(\tilde{\beta},\frac{a}{\beta}\right).\label{basic-ineq-n}
\end {eqnarray}
\end{lemma}
\begin{proof}
First, let us remark that 
\begin{eqnarray*}
\rho(Tx_n,y)&\le&\rho(Tx_n,Ty)+\rho(Ty,y)\le \alpha(x,n,y)\le \beta,\\
\rho(x_n,y)&\le& \alpha(x,n,y)\le \beta,\\
\rho(x_n,Tx_n)&\ge& a= \left(\frac{a}{\beta}\right)\cdot\beta.
\end {eqnarray*}
Moreover,
\begin{eqnarray*}
0<a\le\rho(x_n,Tx_n)&\le & \rho(x_n,y)+\rho(y,Ty)+\rho(Ty,Tx_n)\\
&\le &  2\rho(x_n,y)+\rho(y,Ty)\le 2\alpha(x,n,y)\le 2\beta,
\end {eqnarray*}
so $\displaystyle \frac{a}{\beta}\in(0,2]$.
Thus, we can apply (\ref{uc-ineq-Groetsch}) to obtain
\begin{eqnarray*}
\rho(x_{n+1},y)&=&\rho((1-\lambda_n)x_n\oplus\lambda_nTx_n,y)\\
&\le& \left(1- 2\lambda_n(1-\lambda_n)\eta\left(\alpha(x,n,y),\frac{a}{\beta}\right)\right)\alpha(x,n,y)\\
&=& \rho(x_n,y)+\rho(y,Ty)\\
&&- 2\lambda_n(1-\lambda_n)\eta\left(\alpha(x,n,y),\frac{a}{\beta}\right)\alpha(x,n,y).
\end {eqnarray*} 
Since $\alpha(x,n,y)\le \tilde{\beta}$, and $\eta$ decreases with $r$, we get that 
\[\eta\left(\alpha(x,n,y),\frac{a}{\beta}\right)\ge \eta\left(\tilde{\beta},\frac{a}{\beta}\right),\]
hence
\begin{eqnarray*}
\rho(x_{n+1},y)&\le & \rho(x_n,y)+\rho(y,Ty)\\
&&- 2\lambda_n(1-\lambda_n)\eta\left(\tilde{\beta},\frac{a}{\beta}\right)\alpha(x,n,y)\\
&\le & \rho(x_n,y)+\rho(y,Ty)\\
&&- 2\gamma\lambda_n(1-\lambda_n)\eta\left(\tilde{\beta},\frac{a}{\beta}\right),
\end {eqnarray*}
since $\alpha(x,n,y)\ge \gamma$ by the hypothesis.
\qed\end{proof}

\begin{corollary}\label{cor-eta-tilde}
Assume that $\eta$ decreases with $r$ (for a fixed $\varepsilon$), and moreover that $\eta$ can be written as $\eta(r,\varepsilon)=\varepsilon\cdot\tilde{\eta}(r,\varepsilon)$ such that $\tilde{\eta}$ increases with $\varepsilon$ (for a fixed $r$).\\ 
Let $x,y\in X$, $n\in\N$, and  $\delta,a>0$ be such that 
\begin{equation}
\alpha(x,n,y)\le \delta, \quad  \text{and~} a\le\rho(x_n,Tx_n).\label{hyp-cor-eta-tilde}
\end{equation}
Then
\begin{eqnarray}
\rho(x_{n+1},y)&\le& \rho(x_n,y)+\rho(y,Ty)-2a\lambda_n(1-\lambda_n)\tilde{\eta}\left(\delta,\frac{a}{\delta}\right).
\end {eqnarray}
\end{corollary}
\begin{proof}
Applying Lemma \ref{main-lemma} with $\gamma:=\beta:=\alpha(x,n,y),\tilde{\beta}:=\delta$, we get that 
\begin{eqnarray*}
\rho(x_{n+1},y)&\le& \rho(x_n,y)+\rho(y,Ty)- 2\alpha(x,n,y)\lambda_n(1-\lambda_n)\eta\left(\delta,\frac{a}{\alpha(x,n,y)}\right)\\
&=& \rho(x_n,y)+\rho(y,Ty)- 2a\lambda_n(1-\lambda_n)\tilde{\eta}\left(\delta,\frac{a}{\alpha(x,n,y)}\right)\\
&\le& \rho(x_n,y)+\rho(y,Ty)- 2a\lambda_n(1-\lambda_n)\tilde{\eta}\left(\delta,\frac{a}{\delta}\right),
\end {eqnarray*}
since $\displaystyle \frac{a}\delta\le \frac{a}{\alpha(x,n,y)}$, and $\tilde{\eta}$ increases with $\varepsilon$.
\qed\end{proof}

\begin{corollary}\label{cor-used-Groetsch}
Assume that $\eta$ decreases with $r$ (for a fixed $\varepsilon$).\\ 
Let $x,y\in X, N\in\N$, and $b,c,\gamma,a>0$ be such that
\begin{equation}
\rho(x,y)\le b, \quad  \text{and~}\rho(y,Ty)\le c,
\end{equation}
and for all $n=\overline{0,N}$,
\begin{equation}
\gamma\le \alpha(x,n,y), \quad  \text{and~} a\le\rho(x_n,Tx_n).
\end{equation}
Let $d\ge (N+1)c$.
Then 
\begin{eqnarray}
\rho(x_{N+1},y)&\le& b+(N+1)c-2\gamma\eta\left(b+d,\frac{a}{b+d}\right)\sum_{i=0}^N\lambda_n(1-\lambda_n).
\end {eqnarray}
\end{corollary}
\begin{proof}
Using (\ref{KM-useful-2}), we get for $n=\overline{0,N}$,
\begin{eqnarray*}
\alpha(x,n,y)&=&\rho(x_n,y)+\rho(y,Ty)\le \rho(x,y)+\left(\sum_{i=0}^{n-1}\lambda_i\right)\rho(y,Ty) +\rho(y,Ty)\\
&\le & b+(n+1)\rho(y,Ty)\le b+(N+1)c\le b+d.
\end {eqnarray*}
Applying  Lemma \ref{main-lemma} with $\beta:=\tilde{\beta}:=b+d$, we get 
\begin{eqnarray}
\rho(x_{n+1},y)&\le& \rho(x_n,y)+\rho(y,Ty)- 2\gamma\lambda_n(1-\lambda_n)\eta\left(b+d,\frac{a}{b+d}\right).\label{ineq-cor-1}
\end {eqnarray}
Adding (\ref{ineq-cor-1}) for $n=\overline{0,N}$, it follows that 
\begin{eqnarray*}
\rho(x_{N+1},y)&\le& \rho(x,y)+(N+1)\rho(y,Ty)-2\gamma\eta\left(b+d,\frac{a}{b+d}\right)\sum_{i=0}^N\lambda_n(1-\lambda_n)\\
&\le& b+(N+1)c-2\gamma\eta\left(b+d,\frac{a}{b+d}\right)\sum_{i=0}^N\lambda_n(1-\lambda_n).
\end {eqnarray*}
\qed\end{proof}

\begin{corollary}\label{cor-used-Groetsch-tilde-eta}
In the hypothesis of Corollary \ref{cor-used-Groetsch}, assume moreover that $\eta(r,\varepsilon)$ can be written as $\eta(r,\varepsilon)=\varepsilon\cdot\tilde{\eta}(r,\varepsilon)$ such that $\tilde{\eta}$ increases with $\varepsilon$ (for a fixed $r$).
Then
\begin{eqnarray}
\rho(x_{N+1},y)&\le& b+(N+1)c-2a\tilde{\eta}\left(b+d,\frac{a}{b+d}\right)\sum_{i=0}^N\lambda_n(1-\lambda_n).
\end {eqnarray}
\end{corollary}
\begin{proof}
Follow the proof of Corollary \ref{cor-used-Groetsch}, using this time Corollary \ref{cor-eta-tilde} instead of  Lemma \ref{main-lemma}.
\qed\end{proof}

\section{Main theorem}

\begin{theorem}  \label{main-Groetsch-thm} 
Let $(X,\rho,W)$ be a uniformly convex hyperbolic space with modulus of uniform convexity $\eta$ such that $\eta$ decreases with  $r$(for a fixed $\varepsilon$).\\
Let $C\subseteq X$ be a nonempty convex subset, and $T:C\rightarrow C$ nonexpansive.\\
Assume that $(\lambda_n)$ is a sequence in $[0,1]$ and $\theta :\N\to\N$ is such that for all $n\in\N$, 
\begin{equation}
\sum\limits_{k=0}^{\theta (n)} \lambda_k(1-\lambda_k) \ge n. \label{main-hyp-lambda-gamma}
\end{equation} 
Let $x\in C,b>0$ be such that for any $\delta >0$ there is $y\in C$ with 
\begin{equation}
\rho(x,y)\le b, \quad \mbox{and~~}  \rho(y,Ty) <\delta. \label{main-hyp-x-y}
\end{equation}
Then    
\[\displaystyle\lim_{n\to\infty} \rho(x_n,Tx_n)=0,\]
and moreover
\begin{equation} 
\forall \varepsilon >0\, \forall n\ge h(\varepsilon,b,\theta)\, \big( \rho(x_n,Tx_n)\leq \varepsilon\big), \label{main-thm-conclusion}
\end{equation}
where
\[h(\varepsilon,b,\theta):=\left\{ \begin{array}{ll}
        \displaystyle \theta\left(\left\lceil\frac{b+1}{\varepsilon\cdot\eta\left(b+1,\frac{\varepsilon}{b+1}\right)}
\right\rceil\right) & \text{for~ } \varepsilon <2b\\
        0 & \text{otherwise. }
        \end{array}\right.
\]
\end{theorem}
\begin{proof}
First, let us remark that for any $n\in\N$, and $\delta>0$, by the fact that $(\rho(x_n,Tx_n))$ is non-increasing, and (\ref{main-hyp-x-y})we get
\[\rho(x_n,Tx_n)\le\rho(x,Tx)\le 2\rho(x,y)+\rho(y,Ty)\le 2b+\delta.\]
It follows that $\rho(x_n,Tx_n)\le 2b$ for any $n\in\N$, hence the case $\varepsilon\ge 2b$ is obvious. Let us consider $\varepsilon<2b$ and denote
\begin{eqnarray*}
K&:=&\left\lceil\frac{b+1}{\varepsilon\cdot \eta(b+1,\frac{\varepsilon}{b+1})}
\right\rceil,\\
N&:=&\theta(K)=h(\varepsilon,b,\theta).
\end {eqnarray*}
Assume that for all $n\le N$ we have that $\rho(x_n,Tx_n)>\varepsilon$. Let $\delta>0$ be such that 
\[\delta<\frac{1}{2(N+1)}, \text{~so~}(N+1)\delta<\frac12<1.\]
Let $y\in C$ satisfying (\ref{main-hyp-x-y}) for this $\delta$.
It follows that for all $n=\overline{0,N}$, 
\[\alpha(x,n,y)\ge \frac{\rho(x_n,Tx_n)}2>\frac\varepsilon 2.\]
Applying Corollary \ref{cor-used-Groetsch} with $\displaystyle a:=\varepsilon,\, c:=\delta,\, \gamma:=\frac\varepsilon 2,\, d:=1$, we get 
\begin{eqnarray*}
\rho(x_{N+1},y)&\le& b+(N+1)\delta-\varepsilon\eta\left(b+1,\frac{\varepsilon}{b+1}\right)\sum_{i=0}^N\lambda_n(1-\lambda_n)\\
&\le& b+\frac12-\varepsilon\eta\left(b+1,\frac{\varepsilon}{b+1}\right)\sum_{i=0}^{\theta(K)}\lambda_n(1-\lambda_n)\\
&\le& b+\frac12-\varepsilon\eta\left(b+1,\frac{\varepsilon}{b+1}\right)K, \quad \text{by (\ref{main-hyp-lambda-gamma})}\\
&\le& b+\frac12-(b+1)<0,
\end {eqnarray*}
that is a contradiction.
Thus, there is $n\le h(\varepsilon,b,\theta)$ such that $\rho(x_n,Tx_n)\leq \varepsilon$. Since $(\rho(x_n,Tx_n))$ is non-increasing, (\ref{main-thm-conclusion}) follows.
\qed\end{proof}

\begin{fact}
If $(X,\rho,W)$ is uniformly convex and $\eta$ is any modulus of uniform convexity, then by defining
\begin{equation}
\eta^+(r,\varepsilon):=\inf\{\eta(s,\varepsilon)\mid s\leq r\}, 
\end{equation}
we  get a modulus of  uniform convexity $\eta^+$ which  decreases with  $r$(for a fixed $\varepsilon$), that is, satisfying the requirement from the hypothesis of the above theorem.  
\end{fact}

\begin{fact}\label{tilde-eta}
In the hypothesis of the above theorem, assume moreover that  $\eta(r,\varepsilon)$ can be written as $\eta(r,\varepsilon)=\varepsilon\cdot\tilde{\eta}(r,\varepsilon)$ such that $\tilde{\eta}$ increases with $\varepsilon$ (for a fixed $r$).
Then the bound $h(\varepsilon,b,\theta)$ can be replaced for $\varepsilon <2b$ by
\[ \tilde{h}(\varepsilon,b,\theta):= \theta\left(\left\lceil
\frac{b+1}{2\varepsilon\cdot \tilde{\eta}\left(b+1,\frac{\varepsilon}{b+1}\right)}\right\rceil\right). \] 
\end{fact}
\begin{proof}
Define 
\[K:=\left\lceil\frac{b+1}{2\varepsilon\cdot \tilde{\eta}\left(b+1,\frac{\varepsilon}{b+1}\right)}\right\rceil,\]
and follow the proof of the theorem using Corollary \ref{cor-used-Groetsch-tilde-eta} instead of Corollary \ref{cor-used-Groetsch}.
\qed\end{proof}

As an immediate consequence of our main theorem, we obtain Groetsch's theorem for hyperbolic spaces.

\begin{corollary}
Let $(X,\rho,W)$ be a uniformly convex hyperbolic space, $C\subseteq X$ be a nonempty convex subset, and $T:C\rightarrow C$ nonexpansive such that $T$ has at least one fixed point. 
Assume that $(\lambda_n)$ is a sequence in $[0,1]$ such that 
\[\sum_{n=0}^\infty\lambda_n(1-\lambda_n)=\infty.\]
Then for any $x\in C$,
\[\displaystyle\lim_{n\to\infty} \rho(x_n,Tx_n)=0.\]
\end{corollary}
\begin{proof}
Since $\displaystyle \sum_{n=0}^\infty\lambda_n(1-\lambda_n)$ diverges, for any $n\in\N$ there is $N\in\N$ such that $\displaystyle\sum_{k=0}^N\lambda_k(1-\lambda_k)\ge n,$
so by defining $\theta(n)$ as the least $N$ with this property, (\ref{main-hyp-lambda-gamma}) is satisfied. Let $p$ be a fixed point of $T$.  Then for any $x\in C$, if we take $b:=\rho(x,p)$, (\ref{main-hyp-x-y}) is satisfied with $y:=p$. Hence, we can apply Theorem \ref{main-Groetsch-thm} to get $\displaystyle\lim_{n\to\infty}\rho(x_n,Tx_n)=0$.
\qed\end{proof}

\begin{corollary}\label{bounded-C-general-lambda}
 Let $(X,\rho,W), C,T,(\lambda_n),\theta$ be as in the hypothesis of Theorem \ref{main-Groetsch-thm}. Assume moreover that $C$ is bounded with finite diameter $d_C$.
Then $T$ is $\lambda_n$-asymptotic regular, and  the following holds for all $x\in C$:
\[\forall \varepsilon >0\,\forall n\ge h(\varepsilon,d_C,K)\,
\big(\rho(x_n,Tx_n) <\varepsilon\big), \]
where $h(\eps,d_c,\theta)$ is defined as in Theorem \ref{main-Groetsch-thm} by replacing $b$ with $d_C$.
\end{corollary}
\begin{proof}
If $C$ is bounded, then $C$ has approximate fixed point property, as a consequence of \cite[Theorem 1]{Goebel/Kirk(83)}, which was proved for spaces of hyperbolic type, but, as we have already remarked, any hyperbolic space in our sense is a space of hyperbolic type. It follows that the condition (\ref{main-hyp-x-y}) holds for all $x\in C$ with $d_C$ instead of $b$. Hence, we can apply Theorem \ref{main-Groetsch-thm} for all $x\in C$.
\end{proof}

Thus, for bounded $C$, we get asymptotic regularity for general $(\lambda_n)$, and an explicit bound $h(\eps,d_c,\theta)$ on the rate of asymptotic regularity, which depends only on the error $\varepsilon$, on the diameter $d_C$ of $C$, and on $(\lambda_n)$ only via $\theta$, but not  on the nonexpansive mapping $T$, the starting point $x\in C$ of the iteration or other data related with $C$ and $X$.

The bound $h(\eps,d_C,\theta)$ on the rate of asymptotic regularity can be further simplified in the case of constant $\lambda_n:=\lambda\in(0,1)$.

\bcor\label{bounded-C-constant-lambda}
Let $(X,\rho,W), C, d_C,T$ be as in the hypothesis of Corollary \ref{bounded-C-general-lambda}.
Assume moreover that $\lambda_n:=\lambda\in(0,1)$ for all $n\in\N$.\\
Then $T$ is $\lambda$-asymptotic Al regular, and  for all $x\in C$,
\beq 
\forall \varepsilon >0\,\forall n\ge h(\varepsilon,d_C,\lambda)\,\big(\rho(x_n,Tx_n)< \varepsilon\big), \label{thm-conclusion-constant-lambda}
\eeq
where
\[h(\varepsilon,d_C,\lambda):=\left\{ \ba{ll}\ds\frac{1}{\lambda(1-\lambda)}\left\lceil\frac{d_C+1}{\eps\cdot\eta\left(d_C+1,\ds\frac{\eps}{d_C+1}\right)}\right\rceil & \text{for~ } \varepsilon <2d_C\\
        0 & \text{otherwise. }
        \ea\right.
\]
Moreover, if $\eta(r,\varepsilon)$ can be written as $\eta(r,\varepsilon)=\varepsilon\cdot\tilde{\eta}(r,\varepsilon)$ such that $\tilde{\eta}$ increases with $\eps$ (for fixed $r$), then the bound $h(\varepsilon,d_C,\lambda)$ can be replaced for $\eps<2d_C$ with
\[\tilde{h}(\varepsilon,d_C,\lambda)=\frac{1}{\lambda(1-\lambda)}\left\lceil\frac{d_C+1}{\eps\tilde{\eta}\left(d_C+1,\ds\frac{\eps}{d_C+1}\right)}\right\rceil\]
\ecor
\begin{proof}
(\ref{main-hyp-lambda-gamma}) is satisfied with 
\[\displaystyle \theta:\N\to\N, \, \theta(n):=\frac{n}{\lambda(1-\lambda)}.\]
Hence, we can apply Corollary \ref{bounded-C-general-lambda}, and Remark \ref{tilde-eta}. In this case,

\begin{eqnarray*}
h(\varepsilon,d_C,\theta)&=& \theta\left(\left\lceil\frac{d_C+1}{\eps\cdot\eta\left(d_C+1,\ds\frac{\eps}{d_C+1}\right)}\right\rceil\right)\\
&=&\frac{1}{\lambda(1-\lambda)}\left\lceil\frac{d_C+1}{\eps\cdot\eta\left(d_C+1,\ds\frac{\eps}{d_C+1}\right)}\right\rceil.
\end{eqnarray*}
\end{proof}

As we have proved in Section \ref{UC-hyp}, CAT(0)-spaces are uniformly convex hyperbolic spaces with a ``nice'' modulus of uniform convexity $\ds \eta(r,\varepsilon):=\frac{\varepsilon^2}{8}$, which has the form required in Remark \ref{tilde-eta}. Thus, the above results  can be applied to CAT(0)-spaces. 

\begin{corollary}\label{CAT0-general-lambda}
Let $X$ be a CAT(0)-space,  and $C,d_C,T,(\lambda_n), \theta$ be as in the hypothesis of Corollary \ref{main-Groetsch-thm}. \\
Then $T$ is $\lambda_n$-asymptotic regular, and for all $x\in C$,
\beq 
\forall \varepsilon >0\, \forall n\ge g(\varepsilon,d_C,\theta)\, \big( \rho(x_n,Tx_n)< \varepsilon\big),
\eeq
where
\[g(\varepsilon,d_C,\theta):=\left\{ \ba{ll}
        \ds \theta\left(\left\lceil\frac{8(d_C+1)^2}{\eps^2}\right\rceil\right) & \text{for~ } \varepsilon <2d_C\\
        0 & \text{otherwise. }
        \ea\right.
\]
\end{corollary}
\begin{proof}
Apply Corollary \ref{CAT0-uc}, Corollary \ref{bounded-C-general-lambda}, and Remark \ref{tilde-eta}.
\end{proof}

\begin{corollary}\label{CAT0-constant-lambda}
Let $X,C d_C,T$ be as in the hypothesis of Corollary \ref{CAT0-general-lambda}. 
Assume that $(\lambda_n)=\lambda\in(0,1)$.
Then $T$ is $\lambda$-asymptotic regular, and for all $x\in C$,
\beq 
\forall \varepsilon >0\, \forall n\ge \tilde{g}(\varepsilon,d_C,\lambda)\, \big( \rho(x_n,Tx_n)< \varepsilon\big),
\eeq
where
\[\tilde{g}(\varepsilon,d_C,\lambda):=\left\{ \ba{ll}
        \ds\frac{1}{\lambda(1-\lambda)}\left\lceil\frac{8(d_C+1)^2}{\eps^2}\right\rceil & \text{for~ } \varepsilon <2d_C\\
        0 & \text{otherwise. }
        \ea\right.
\]
\end{corollary}
\begin{proof}
Apply Corollary \ref{CAT0-uc}, Corollary \ref{bounded-C-constant-lambda}, and Remark \ref{tilde-eta}.
\end{proof}

Hence, in the case of CAT(0)-spaces, we get a quadratic (in $1/varepsilon$) rate of asymptotic regularity.\\[0.1cm]

\noindent {\bf Acknowledgements}\\

\noindent I would like to express my gratitude to Ulrich Kohlenbach for bringing the problem under consideration to my attention and for the numerous discussions we had on the subject.

\end{document}